\documentclass{amsart}
\newtheorem{theorem}{Theorem}[section]
\newtheorem{lemma}[theorem]{Lemma}
\theoremstyle{definition}
\newtheorem{definition}[theorem]{Definition}
\newtheorem{example}[theorem]{Example}
\newtheorem{corollary}[theorem]{Corollary}
\newtheorem{proposition}[theorem]{Proposition}

\theoremstyle{remark}
\newtheorem{remark}[theorem]{Remark}
\numberwithin{equation}{section}

\newcommand{\eps}{\varepsilon}

\usepackage{color,graphicx,amssymb}
\usepackage{graphicx}
\begin{document}
\title{SRB Measures for Certain Markov Processes}
\author{Wael Bahsoun}
    \address{Department of Mathematical Sciences, Loughborough University, 
Loughborough, Leicestershire, LE11 3TU, UK}
\email{W.Bahsoun@lboro.ac.uk}
\author{Pawe\l\ G\'ora}
\address{Department of Mathematics and Statistics, Concordia University, Montreal, Quebec, H3G 1M8, Canada}
\email{pgora@mathstat.concordia.ca}
\subjclass{Primary 37A05, 37E05, 37H99}
\keywords{Iterated Function System, SRB-Measures}
\begin{abstract}
We study Markov processes generated by iterated function systems (IFS). The constituent maps of the IFS are monotonic transformations of the interval. We first obtain an upper bound on the number of SRB (Sinai-Ruelle-Bowen) measures for the IFS. Then, when all the constituent maps have common fixed points at 0 and 1, theorems are given to analyze properties of the ergodic invariant measures $\delta_0$ and $\delta_1$. In particular, sufficient conditions for $\delta_0$ and/or $\delta_1$ to be, or not to be, SRB measures are given. We apply some of our results to asset market games.
\end{abstract}
\maketitle
\section{Introduction}
In the 1970's, Sinai, Ruelle and Bowen studied the existence of an important class of invariant measures in the context of deterministic dynamical systems. These invariant measures are nowadays known as SRB (Sinai-Ruelle-Bowen) measures \cite{Y}. SRB measures are distinguished among other ergodic invariant measures because of their physical importance. In fact, from ergodic theory point of view, they are the only useful ergodic measures. This is due to the fact that SRB measures are the only ergodic measures for which the Birkhoff  Ergodic Theorem holds on a set of positive measure of the phase space.  In this note, we study SRB measures in a stochastic setting---Markov processes generated by iterated function systems (IFS). 

\bigskip

An IFS\footnote{In some of the literature an IFS is called a \textit{random map} or a \textit{random transformation}.} is a discrete-time random dynamical system \cite{Ar, Ki} which consists of a finite collection of transformations and a probability vector $\{\tau_s;p_s\}_{s=1}^{L}$. At each time step, a transformation $\tau_s$ is selected with probability $p_s>0$ and applied to the process. IFS has been a very active topic of research due to its wide applications in fractals and in learning models. The survey articles \cite{Di, St} contain a considerable list of references and results in this area. 

\bigskip

The systems which we study in this note do not fall in the category of the IFS\footnote{ In most articles about IFS, the constituent maps are assumed to be contracting or at least contracting on average. Here we do not impose any assumption of this type. In fact the class of IFS which we study in Section \ref{martingale} cannot satisfy such assumptions.} considered in \cite{Di, St} and references therein. Moreover, in general, our IFS do not satisfy the classical splitting\footnote{In particular, when all the maps have common fixed points at 0 and 1. See Section \ref{martingale}.} condition of \cite{DF}. In fact, our aim in this note is to depart from the traditional goal of finding sufficient conditions for an IFS to admit a \textit{unique attracting invariant measure} \cite{DF,Di,St}. Instead, we study cases where an IFS may admit more than one invariant measure and aim to identify the \textit{physically relevant} ones; i.e., invariant measures for which the Ergodic Theorem holds on a set of positive measure of the ambient phase space. We call such invariant measures SRB.

\bigskip

Physical SRB measures for random maps have been studied by Buzzi \cite{Bu} in the context of random Lasota-Yorke maps. However, Buzzi's definition of a basin of an SRB measure is different from ours. We will clarify this difference in Section \ref{Pre}. A general concept of an SRB measure for general random dynamical systems can be found in the survey article \cite{Liu}. In this note we study physical SRB measures for IFS whose constituent maps are strictly increasing transformations of the interval. We obtain an upper bound on the number of SRB measures for the IFS. Moreover, when all the constituent maps have common fixed points at 0 and 1, we provide sufficient conditions for $\delta_0$ and/or $\delta_1$ to be, or not to be, SRB measures. To complement our theoretical results, we show at the end of this note that examples of IFS of this type can describe evolutionary models of financial markets \cite{EHS}. 

\bigskip

In Section \ref{Pre}  we introduce our notation and main definitions. 
In particular, Section \ref{Pre} includes the definition of an SRB measure for an IFS. 
In Section \ref{SRB:bound} we identify the structure of the basins of SRB measures and we obtain a sharp upper bound on the number of SRB measures. Section \ref{martingale} 
contains sufficient conditions for $\delta_0$ and $\delta_1$,  the delta measures concentrated at $0$ and $1$ respectively, to be SRB. It also contains sufficient conditions for $\delta_0$ and $\delta_1$ not to be SRB measures. Our main results in this section are Theorems \ref{main} and Theorem \ref{NotSRB}. 
In Section \ref{information} we study ergodic properties of $\delta_0$ and $\delta_1$ without having any information about the probability vector of the IFS. 
In Section \ref{games} we apply our results to asset market games.  In particular, we find a generalization of the famous Kelly rule \cite{Ke} which expresses the principle of ``betting your beliefs''. The importance of our generalization lies in the fact that it does not require the full knowledge of the probability distribution of the random states of the system. Section \ref{app} contains an auxiliary result which we use in the proof of Theorem \ref{NotSRB}.
 
\section{Preliminaries}\label{Pre}
\subsection{Notation and assumptions} Let $([0,1],\mathfrak B)$ be the measure space where $\mathfrak B$ is the Borel
$\sigma$-algebra on $[0,1]$. Let $\lambda$ denote Lebesgue measure on $([0,1],\mathfrak B)$ and $\delta_r$ denote the delta measure concentrated at point $r\in [0,1]$. Let $S=\{1,\dots ,L\}$ be a finite set and $\tau_s$, $s\in S$, be continuous transformations from the unit interval into itself. We assume:\\

\noindent (A) $\tau_s$ are strictly increasing.\\

Let $\mathbb {\bf p}=(p_s)_{s=1}^{L}$ be a probability vector on $S$ such that for all $s\in S$, $p_s>0$.
The collection 
$$F=\{\tau_1,\tau_2,\dots ,\tau_L;p_1,p_2,\dots ,p_L\}$$ 
is called  an iterated function system (IFS) with probabilities. 

\bigskip

We denote the space of sequences $\omega=\{s_1,s_2,\dots\}$, $s_l\in S$, by $\Omega$. The topology on $\Omega$ is defined as the product of the discrete topologies on $S$. Let $\pi_{\bf p}$ denote the Borel measure on $\Omega$ defined as the product measure ${\bf p}^{\mathbb N}$. Moreover, we write 
$$s^t:=(s_1,s_2,\dots,s_t)$$ 
for the history up to time $t$, and for any $r_0\in[0,1]$ we write
$$r_t(s^t):=\tau_{s_t}\circ \tau_{s_{t-1}}\circ\cdots\circ\tau_{s_1}(r_0).$$ 
Finally, by $E(\cdot)$ we denote the expectation with respect to $\bf{p}$, by $E(\cdot|s^t)$ the conditional expectation given the history up to time $t$ and by ${\rm var}(\cdot)$ the variance with respect to $\bf{p}$.
\subsection{Invariant measures}
$F$ is understood as a Markov process with a transition function 
$$\mathbb P (r, A)=\sum_{s=1}^{L}p_s\chi_{A}(\tau_s(r)),$$
where $A\in\mathfrak B$ and $\chi_{A}$ is the characteristic function of the set $A$. The transition function $\mathbb P$ induces an operator $P$ on measures on $([0,1], \mathfrak B)$ defined by
\begin{equation}
\begin{split}
P\mu(A)&=\int_0^1\mathbb P(r,A) d\mu(r)\\
&=\sum_{s=1}^{L}p_{s}\mu(\tau^{-1}_{s}A).\\
\end{split}
\end{equation}
Following the standard notion of an invariant measure for a Markov process, we call a probability measure $\mu$ on $([0,1],\mathfrak B)$ $F$-invariant probability measure if and only if 
$$P\mu=\mu.$$ 
Moreover, it is called ergodic if it cannot be written as a convex combination of other invariant probability measures.
\subsection{SRB measures}
Let $\mu$ be an ergodic probability measure for the IFS. Suppose there exists a set of positive Lebesgue measure in $[0,1]$ such that
\begin{equation}\label{SRBbasin}
\frac{1}{T}\sum_{t=0}^{T-1}\delta_{r_t(s^t)}\overset{\text{weakly}}{\to}\mu\hskip 1cm\text{with } \pi_{\bf p}\text{-probability one.}
\end{equation}
Then $\mu$ is called an SRB (Sinai-Ruelle-Bowen) measure. The set of points $r_0\in[0,1]$ for which (\ref{SRBbasin}) is satisfied will be called the basin\footnote{Our definition of a basin is different from Buzzi's definition \cite{Bu}. In his definition he defines \textit{random} basins $B_{\omega}(\mu)$ for an SRB measure. In particular, according to Buzzi's definition, for the same SRB measure, basins corresponding to two different $\omega$'s may differ on a set of positive lebsegue measure of $[0,1]$. See \cite{Bu} for more detials.} of $\mu$ and it will be denoted by $B(\mu)$. Obviously, if $\lambda(B(\mu))=1$ then $\mu$ is the unique SRB measure of $F$.

\section{Number of SRB measures and their basins}\label{SRB:bound}
The basin of an SRB measure for the systems we are dealing with is described by the following two propositions.
\begin {proposition}\label{Basin}
Let $\mu$ be an SRB measure and $B(\mu)$ be its basin. Let $r_0, \bar r_0\in B(\mu)$, $r_0> \bar r_0 $. Then $ [\bar r_0, r_0]\subseteq B(\mu)$.  
\end{proposition}
\begin{proof} When weak convergence is considered on an interval, then $\mu_n \overset{\text{weakly}}{\to}\mu$
if and only if $\mu_n(f)\to\mu(f)$ for any $C^1$ function\footnote{Here is a sketch of the proof of this claim:
Assume $$\mu_n(f)\to\mu(f)\,$$ for any $f\in C^1([0,1])$. Let $g$ be a continuous function and let 
$\{f_k\}_{k\ge 1}$ be a sequence of $C^1$ functions converging to $g$ in $C^0$ norm.
We have
\begin{equation*}
\begin{split}
|\mu_n(g)-\mu(g)|&\le |\mu_n(g)-\mu_n(f_k)|+|\mu_n(f_k)-\mu(f_k)|+|\mu(f_k)-\mu(g)|\\
&\le 2 \|f_k-g\|_{C^0}+ |\mu_n(f_k)-\mu(f_k)|\ .
\end{split}
\end{equation*}
Now, for any $\varepsilon>0$, we can find $k_0$ such that
$2 \|f_{k_0}-g\|_{C^0}<\varepsilon/2$ and then we can find $n_0$ such that for any $n\ge n_0$ we have
$|\mu_n(f_{k_0})-\mu(f_{k_0})|<\varepsilon/2$.}. 
Since every $C^1$ function is a difference of two continuous increasing functions, 
this means that $\mu_n \overset{\text{weakly}}{\to}\mu$
if and only if $\mu_n(f)\to\mu(f)$ for any continuous increasing function. 

\bigskip

Let $r_0, \bar r_0\in B(\mu)$ and $\bar r_0 <r'_0<r_0$. We will show that $r'_0\in B(\mu)$.
Let assume that $f$ is continuous and increasing. Let us fix an $s^t$ for which
$$\lim_{T\to\infty}\frac{1}{T}\sum_{t=0}^{T-1}f(\bar r_t(s^t))= \lim_{T\to\infty}\frac{1}{T}\sum_{t=0}^{T-1}f(r_t(s^t))=
 \mu(f).$$
We have  $\bar r_t(s^t)<r'_t(s^t)<r_t(s^t)$ (since all $\tau_s$ are increasing) and
$$\frac{1}{T}\sum_{t=0}^{T-1}f(\bar r_t(s^t))\le \frac{1}{T}\sum_{t=0}^{T-1}f(r'_t(s^t))\le
 \frac{1}{T}\sum_{t=0}^{T-1}f(r_t(s^t)).$$
The averages on the left and on the right have common limit $\mu(f)$. Thus, 
$$\frac{1}{T}\sum_{t=0}^{T-1}\delta_{r'_t(s^t)}(f)=\frac{1}{T}\sum_{t=0}^{T-1}f(r'_t(s^t))\to \mu(f).$$
Since the event
$$\{\lim_{T\to\infty}\frac{1}{T}\sum_{t=0}^{T-1}f(\bar r_t(s^t))= \lim_{T\to\infty}\frac{1}{T}\sum_{t=0}^{T-1}f(r_t(s^t))
= \mu(f)\}$$
occurs with $\pi_{\bf p}$-probability 1, the event
$$\{\frac{1}{T}\sum_{t=0}^{T-1}f(r'_t(s^t))\to \mu(f)\}$$
also occurs with $\pi_{\bf p}$-probability 1. 
\end{proof} 

\begin{figure}[h] 
  \centering
 \includegraphics[bb=0 0 400 400,width=4.02in,height=4.02in,keepaspectratio]{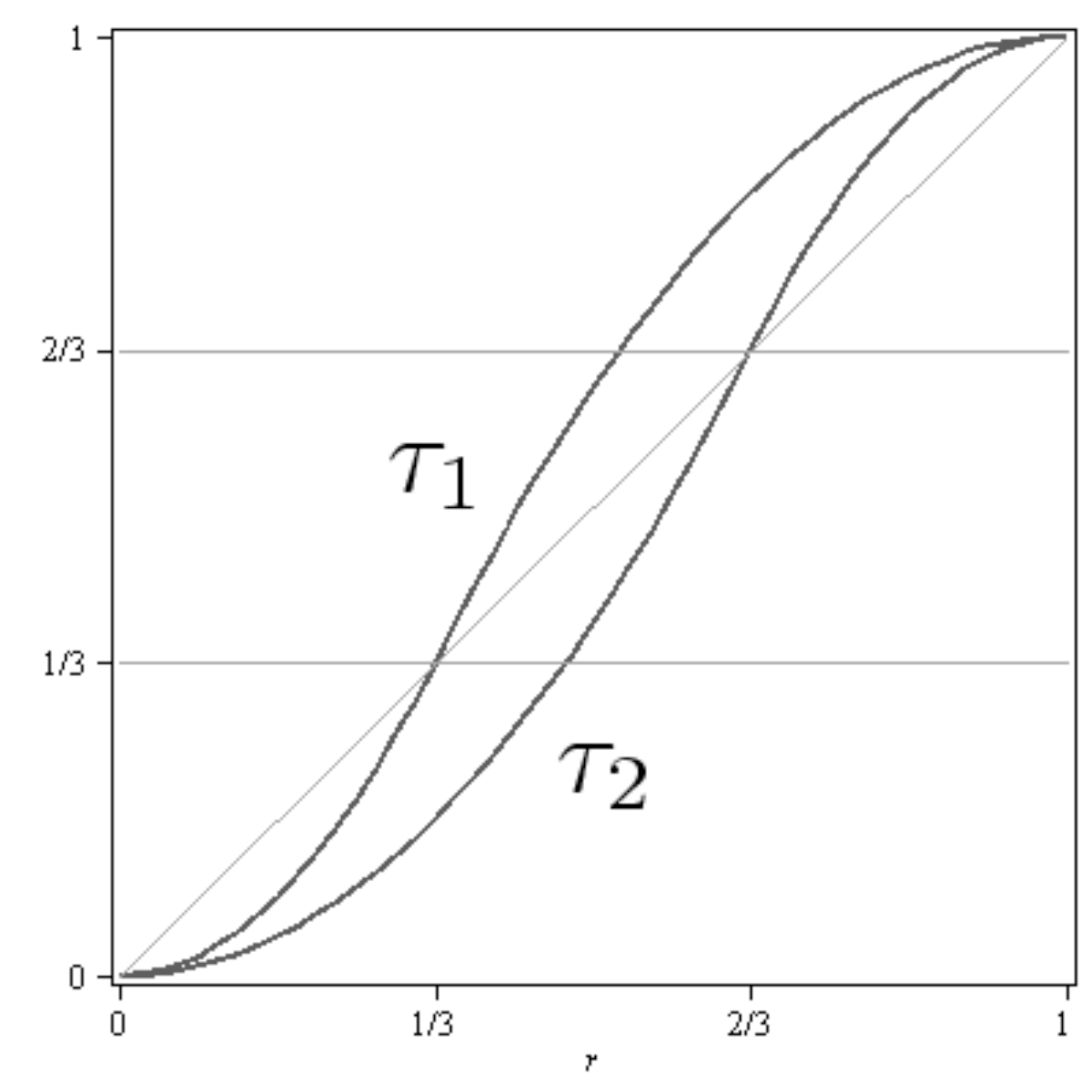}
  \caption{Maps $\tau_1$ and $\tau_2$ in Example \ref{Ex:basin}}
  \label{fig:Example}
\end{figure}

\begin {proposition}\label{Basin_ends}
Let $\mu$ be an SRB measure and $B(\mu)=\langle a,b\rangle$ be its basin, where $\langle a,b\rangle$
denotes an interval closed or open at any of the endpoints. Then, 
$$ \tau_s(a)\ge a \ \ ,  \ \ s=1,\dots,L,  \ \ \text{and if }\\ a\not=0\ \ \text{then }\\ \tau_s(a)= a \ \ \text{for at least one } s;$$
$$ \tau_s(b)\le b \ \ ,  \ \ s=1,\dots,L, \ \ \text{and if }\\ b\not=1\ \ \text{then }\ \ \tau_s(b)= b \ \ \text{for at least one } s.$$
\end{proposition}
\begin{proof} We will prove only the second claim with $b\not=1$. The first claim is proven analogously.

Assume that $\tau_{s_0}(b)> b$, for some $1\le s_0\le L$. Then, we can find $r_0\in(a,b)$ such that $\tau_{s_0}(r_0)> b$.
For arbitrary continuous function $f$, for $\omega\in A\subset \Omega$ with  $\pi_{\bf p}(A)=1$, we have 
$$\lim_{T\to\infty}\frac{1}{T}\sum_{t=0}^{T-1}f(r_t(s^t))= \mu(f).$$
The set $A_{s_0}=\{(s_1,s_2,\dots): (s_0,s_1,s_2,\dots)\in A\}$ is also of $\pi_{\bf p}$-probability 1. Let $r'_0=\tau_{s_0}(r_0)$
and let $(s^t)'$ denote the initial subsequences of length $t$ of $\omega\in A_{s_0}$ .
Then,
$$\frac{1}{T}\sum_{t=0}^{T-1}f(r'_t((s^t)'))= \frac{1}{T}\sum_{t=0}^{T-1}f(r_t(s^t))-\frac{1}{T}f(r_0)+\frac{1}{T}
f(r'_{T-1}((s^{T-1})'))\underset{T\to +\infty}{\longrightarrow} \mu(f).$$
This shows that $\tau_{s_0}(r_0)\in B_\omega(\mu)$ and contradicts the assumptions.

Now, we assume that $ \tau_s(b)< b$,  $s=1,\dots,L $. Then, we can find $r_0>b$ such that $\tau_{s}(r_0)\in (a,b)$ for all $s$.
Let $$A_s=\{\omega: \lim_{T\to\infty}\frac{1}{T}\sum_{t=0}^{T-1}f(r'_t(s^t))= \mu(f),\ \text{for}\ r'_0=\tau_{s}(r_0)\}\ ,\ s=1,\dots,L \ .$$
We have $\pi_{\bf p}(A_s)=1$ for each $s$. Hence, $\pi_{\bf p}(A)=1$, where $A=\cup_{1\le s\le L} (s,A_s)$
and $(s,A_s)=\{(s,s_1,s_2,s_3,\dots): (s_1,s_2,s_3,\dots)\in A_s\}$.
For arbitrary continuous function $f$, for $\omega\in A$, if $\omega_1=s$ we have 
\begin{equation*}\begin{split}
&\lim_{T\to\infty}\frac{1}{T}\sum_{t=0}^{T-1}f(r_t(s^t))\\
&=\lim_{T\to\infty}\left(\frac{1}{T}\sum_{t=0}^{T-1}f(r'_t((s^t)'))
+\frac{1}{T}f(r_0)-\frac{1}{T}f(r'_{T-1}((s^{T-1})'))\right)= \mu(f),
\end{split}
\end{equation*}
where $r'_0=\tau_s(r_0)$ and $(s^t)'$ are the initial subsequences of length $t$ of $\omega\in A_s$. This implies that $r_0\in B(\mu)$. Since $r_0>b$, this leads to a contradiction. 
\end{proof}
We now state the main result of this section. Firstly, we recall that  $\langle \cdot ,\cdot \rangle$  denotes an interval which is closed or open at any of the endpoints. Secondly, we define a 
set \textbf{BS} whose elements are intervals of the form $\langle \cdot ,\cdot \rangle$ with the following property:
$$\langle a ,b \rangle \in \textbf{BS}$$
if and only if 
$$ \tau_s(a)\ge a \ \ ,  \ \ s=1,\dots,L  \ \ \text{and}\ \ \tau_s(a)= a \ \ \text{for at least one } s;$$
and 
$$ \tau_s(b)\le b \ \ ,  \ \ s=1,\dots,L  \ \ \text{and}\ \ \tau_s(b)= b \ \ \text{for at least one } s.$$
\begin{theorem}\label{NumSRB}
The number of SRB measures of $F$ is bounded above by the cardinality of the set \textbf{BS}. In particular, if $0$ and $1$ are the only fixed points of some $\tau_{s_0}$, $s_0\in S$,  then $F$ admits at most one SRB measure.
\end{theorem}
\begin{proof}
The fact that number  of SRB measures of $F$ is bounded above by the cardinality of the set \textbf{BS} is a direct consequence of Proposition \ref{Basin_ends}. To elaborate on the second part of the theorem, assume without loss of generality that $\tau_{s_0}(r)>r$ for all $r\in(0,1)$. Obviously, by Proposition \ref{Basin_ends}, if all the other maps $\tau_s$, $s\in S\setminus\{s_0\}$ has no fixed points in $(0,1)$, then $F$ admits at most one SRB measure. So let us assume that there exists an $s^*\in S\setminus\{s_0\}$ such that $\tau_{s^*}$ has a finite or infinite number of fixed points in $[0,1]$. In the case of finite number of fixed points, denote the fixed points of $\tau_{s^*}$ in [0,1] by $r^*_i$, $i=1,\dots ,q$, such that $0\le r^*_1<r^*_2<\cdots <r^*_q\le 1$. Since $\tau_s(r^*_i)>r^*_i$ for all $r^*_i\in (0,1)$, the only possible basin for an SRB measure would be either $\langle r^*_{q-1},1\rangle $ or $\langle r^*_q,1\rangle $. In the case of infinite number of fixed points, let
$$\bar r=\sup\{r\in(0,1):\, \tau_{s^*}(r)=r\}.$$
If $\bar r<1$, then $\tau_{s_0}(\bar r)>\bar r$. By Proposition \ref{Basin_ends}, $\langle \bar r,1\rangle $ is the only possible basin for an SRB measure. If $\bar r =1$, let $\bar J$ denote the closure of the set of fixed points of $\tau_{s^*}$ and let $\bar J_0\subseteq \bar J$ be the minimal closed subset of $\bar J$ which contains the point $1$. $\bar J_0$ is the only possible basin for an SRB measure. Moreover, it cannot be decomposed into basins of different SRB measures. Indeed, let $J_1\cup J_2=\bar J_0$ such that $J_1= \langle a,b\rangle$ with $b<1$. Since $\tau_{s_0}(b)>b$, by Proposition \ref{Basin_ends}, $J_1$ cannot be a basin of an SRB measure. Thus, $F$ admits at most one SRB measure. 
\end{proof}
The following example shows that Proposition \ref{Basin_ends} can be used to identify intervals which are not in the basin of an SRB measure. In particular, it shows that the bound obtained on the number of SRB measures in Theorem \ref{NumSRB} is really sharp. 
\begin{example}\label{Ex:basin}
Let
\begin{equation*}
\tau_1(r)=  \begin{cases} 3r^2     &,\  \text {for}\ \ 0\le r\le 1/3;\\
                          1-\frac 32(r-1)^2     &,\  \text {for}\ \ 1/3< r\le 1;
            \end{cases} \ \ \ , \ \ \
\end{equation*}
and
\begin{equation*} 
\tau_2(r)=  \begin{cases} \frac 32 r^2    &,\  \text {for}\ \ 0\le r\le 2/3;\\
                          1-3(r-1)^2   &,\  \text {for}\ \ 2/3< r\le 1.
            \end{cases}
\end{equation*}
The graphs of the above maps are shown in Figure \ref{fig:Example}. Using Proposition \ref{Basin_ends}, we see that the points of the interval $(1/3,2/3)$ do not belong to a basin of any SRB measure. Moreover,  by Theorem \ref{NumSRB}, $F$ admits at most two SRB measures. Indeed, one can easily check that $\delta_0$ and $\delta_1$ are the only SRB measures with basins $B(\delta_0)=[0,1/3]$ and $B(\delta_1)=[2/3,1]$ respectively.  For any $r\in [0,1/3)$ for all $\omega$'s the averages $\frac{1}{T}\sum_{t=0}^{T-1}\delta_{r_t(s^t)}$ converge weakly to $\delta_0$. For $r=1/3$ the only $\omega$ for which this does not happen is $\omega=\{1,1,1,\dots\}$ so again the averages converge
weakly to $\delta_0$ with $\pi_{\bf p}$-probability 1. Similarly, we can show that $B(\delta_1)=[2/3,1]$. If $r\in (1/3,2/3)$, then with positive $\pi_{\bf p}$-probability the averages converge to $\delta_0$ and with
positive $\pi_{\bf p}$-probability the averages converge to $\delta_1$. Thus, these points do not belong to a basin of any SRB measure and there are only two SRB measures.
\end{example}
\section{Properties of $\delta_0$ and $\delta_1$}\label{martingale}
In addition to condition (A), we assume in this section that for all $s\in S$:\\

\noindent (B) $\tau_s(0)=0$ and $\tau_s(1)=1$;\\

\bigskip
 
Obviously by Condition (B) the delta measures $\delta_0$ and $\delta_1$ are ergodic probability measures for the IFS. We will be mainly concerned with the following question: When does $F$ have $\delta_0$ and/or $\delta_1$ as SRB measures? We start our analysis by proving a lemma which provides a sufficient condition for $\delta_{x}$, the point measure concentrated at $x\in[0,1]$, to be an SRB measure.
\begin{lemma}\label{fix}
Suppose that $\tau_s(x)=x$ for all $s\in\{1,\dots ,L\}$ and that there exists an initial point of a random orbit $r_0$, $r_0\not= x$,  for which $\lim_{t\to\infty}r_t(s^t)=x$ with probability $\pi_{\bf p}=1$. Then $\delta_x$ is an SRB measure for $F$ and $B_\omega(\delta_x)\supseteq [x,r_0]$\footnote{The notation here is for the case when $r_0>x$.}.
\end{lemma}
\begin{proof}
Let $f$ be a continuous function on $[0,1]$. Let $r_0\not=x$ and fix a history $s^t$ for which $\lim_{t\to\infty}r_t(s^t)=x$. Then
$$\lim_{t\to\infty}f(r_t(s^t))=f(x).$$
Consequently
$$\lim_{T\to\infty}\frac{1}{T}\sum_{t=0}^{T-1}f(r_t(s^t))=f(x).$$
Since the event 
$$\{\lim_{t\to\infty}r_t(s^t)=x\}$$
appears with probability one, the event
$$\{\lim_{T\to\infty}\frac{1}{T}\sum_{t=0}^{T-1}f(r_t(s^t))=f(x)\}$$
also appears with probability one. Thus, by Proposition \ref{Basin}, $\delta_x$ is an SRB measure for $F$ and $B_\omega(\delta_x)\supseteq [x,r_0]$.
\end{proof}
The following lemma, which is easy to prove, is a key observation for our main results in this section. 
\begin{lemma}\label{key}
Each constituent map of the IFS can be represented as follows:
$$\tau_s(r)=r^{\beta_s(r)},$$
with $ \beta_s(r)$ satisfying:
\begin{enumerate}
\item $\beta_s(r)>0$ in $(0,1)$ ;
\item $(\ln r) \beta_s(r)$ increasing;
\item $\lim_{r\to 0}(\ln r)\beta_s(r)=-\infty$; 
\item $\lim_{r\to 1}(\ln r)\beta_s(r)=0.$
\end{enumerate}
\end{lemma}
In the rest of this section, the following notation will be used:
$$\alpha_{t}\overset{\text{def}}{:=}\beta_s(r_{t-1})\text{ with probability }p_s,\, t=1,2,\dots$$
\begin{theorem}\label{main}
Let $F=\{\tau_s;p_s\}_{s\in S}$ be an IFS such that $\tau_s(r)=r^{\beta_s(r)}$. Assume that 
$0<b_s\le\beta_s(r)\le B_s<\infty$ for all $r\in [0,1]$.
\begin{enumerate}
\item If $E(\ln\alpha_{t}|s^{t-1})\le 0$ a.s., then $\lim_{t\to\infty}r_t(s^t)\not=0$ a.s.  
\item If  $\limsup_{T\to\infty}\frac{1}{T}\sum_{t=1}^{T}E(\ln\alpha_{t}|s^{t-1})<0$ a.s., then $\lim_{t\to\infty}r_t(s^t)=1$ a.s.
\item If $\liminf_{T\to\infty}\frac{1}{T}\sum_{t=1}^{T}E(\ln\alpha_{t}|s^{t-1})>0$ a.s., then $\lim_{t\to\infty}r_t(s^t)=0$ a.s.
\end{enumerate} 
\end{theorem}
\begin{proof}
Let us consider the sequence of random exponents
$$\alpha(t)=\alpha_{t}\alpha_{t-1}\cdots\alpha_{2}\alpha_{1},$$
where $\alpha_{i}=\beta_s(r_{i-1})$ with probability $p_s$,
and observe that
$$ r_{t}(s^t)=r^{\alpha(t)}.$$
We have
$$\ln\alpha(t+1)=\ln\alpha_{t+1}+\ln\alpha(t),$$
and, with probability one,
$$E(\ln\alpha(t+1)|s^t)-\ln(\alpha(t))=E(\ln\alpha_{t+1}|s^t)\le 0.$$
Therefore, $\ln\alpha(t)$ is a supermartingale. Moreover, because $0<b_s\le\beta_s(r_t)\le B_s<\infty$, $|\ln\alpha(t+1)-\ln\alpha(t)|=|\ln\alpha_{t+1}|<\infty$. Hence $\ln\alpha(t)$ is a supermartingale with bounded increments. Thus, using Theorem 5.1 in Chapter VII of \cite{Sh}, with probability one $\ln\alpha(t)$ does not converge to $+\infty$. Consequently, with probability one, $ r_{t}(s^t)=r^{\alpha(t)}$ does not converge to zero.

We now prove the second statement of the theorem. Again we consider the sequence of random exponents
$$\alpha(t)=\alpha_{t}\alpha_{t-1}\cdots\alpha_{2}\alpha_{1}.$$
Let $M_t$ denote the martingale difference
$$
M_t:=\ln\alpha_t-E(\ln\alpha_{t}|s^{t-1}).
$$
We have $E(M_t)=0$ and $\ln\alpha_t$ is uniformly bounded. Therefore, by the strong law of large numbers (see Theorem 2.19 in \cite{HH}), with probability one
\begin{equation}\label{eq:MD}
\lim_{T\to\infty}\frac{1}{T}\sum_{t=1}^{T}M_t=0.
\end{equation}
Therefore, with probability one,
\begin{equation*}
\begin{split}
\limsup_{T\to\infty}\frac{1}{T}\ln\alpha(T)&=\limsup_{T\to\infty}\frac{1}{T}\sum_{t=1}^{T}\ln\alpha_{t}\\
&=\limsup_{T\to\infty}\frac{1}{T}\sum_{t=1}^{T}M_t +\limsup_{T\to\infty}\frac{1}{T}\sum_{t=1}^{T}E(\ln\alpha_{t}|s^{t-1})<0.
\end{split}
\end{equation*}
From this we can conclude that for $T$ large enough there is a positive random variable $\eta$ such that
$$\alpha(T)\le e^{-T\eta} \text{ a.s.}$$
Thus, since $r\in [0,1]$, for $T$ large enough we obtain
$$r_{T+1}=r^{\alpha(T)}\ge r^{e^{-T\eta}} \text{ a.s.}$$
By taking the limit of $T$ to infinity we obtain
$$\lim_{T\to\infty}r_{T+1}=\lim_{T\to\infty}r^{\alpha(T)}\ge\lim_{T\to\infty} r^{e^{-T\eta}}=1 \text{ a.s.}$$
The proof of the third statement is very similar to the proof of the second one with slight changes. In particular,  using (\ref{eq:MD}), we see that, with probability one, 
$$\liminf_{T\to\infty}\frac{1}{T}\ln\alpha(T)>0.$$
From this we can conclude that for $T$ large enough there is a positive random variable $\eta$ such that
$$\alpha(T)\ge e^{T\eta} \text{ a.s.}$$
Thus, since $r\in [0,1]$, for $T$ large enough we obtain
$$r_{T+1}=r^{\alpha(T)}\le r^{e^{T\eta}} \text{ a.s.}$$
By taking the limit of $T$ to infinity we obtain
$$\lim_{T\to\infty}r_{T+1}=\lim_{T\to\infty}r^{\alpha(T)}\le\lim_{T\to\infty} r^{e^{T\eta}}=0 \text{ a.s.}$$
\end{proof}
\begin{corollary}\label{Cor1}
Let $F=\{\tau_s;p_s\}_{s\in S}$ be an IFS such that $\tau_s(r)=r^{\beta_s(r)}$. Assume that $0<b_s\le\beta_s(r)\le B_s<\infty$ for all $r\in [0,1]$.
\begin{enumerate}  
\item If   $\limsup_{T\to\infty}\frac{1}{T}\sum_{t=1}^{T}E(\ln\alpha_{t}|s^{t-1})<0$ a.s., then $\delta_1$ is the unique SRB measure of $F$ with $B(\delta_1)=(0,1]$ .
\item If  $\liminf_{T\to\infty}\frac{1}{T}\sum_{t=1}^{T}E(\ln\alpha_{t}|s^{t-1})>0$ a.s., then $\delta_0$ is the unique SRB measure of $F$ with $B(\delta_0)=[0,1)$.
\end{enumerate} 
\end{corollary}
\begin{proof}
The proof is a consequence of statements (2) and (3) of Theorem \ref{main} and Lemma \ref{fix}.
\end{proof}
\begin{remark}
Observe that:
\begin{enumerate}
\item $\sum_sp_s\ln B_s\le 0\implies E(\ln\alpha_{t}|s^{t-1})\le 0$ a.s.
\item $\sum_sp_s\ln B_s< 0\implies \limsup_{T\to\infty}\frac{1}{T}\sum_{t=1}^{T}E(\ln\alpha_{t}|s^{t-1})<0$ a.s.
\item $\sum_sp_s\ln b_s> 0\implies \liminf_{T\to\infty}\frac{1}{T}\sum_{t=1}^{T}E(\ln\alpha_{t}|s^{t-1})>0$ a.s. 
\end{enumerate}
Thus, the conditions in the statements of Theorem \ref{main} and Corollary \ref{Cor1}  are very easy to check for certain systems.
\end{remark}
\begin{remark}
In the proof of statement (1) of Theorem \ref{main}, we have with probability $\pi_{\bf p}=1$, $\lim_{t\to\infty}\ln\alpha(t)\not=\infty$. In general, it is not clear that this implies that $\delta_0$ is not an SRB measure. However, in the following theorem under additional natural assumption on the variance of $\ln\alpha_t$ we show that $\delta_0$ is indeed not an SRB measure. 
\end{remark}
\begin{theorem}\label{NotSRB}
If $E(\ln\alpha_t|s^{t-1})\le 0$ and ${\rm var}(\ln\alpha_t|s^{t-1})\ge d > 0$, for all $t\ge 1$, then $\delta_0$ is not an SRB measure of $F$.
\end{theorem}
\begin{proof}
Consider the sequence of random exponents
$$\alpha(t)=\alpha_{t}\alpha_{t-1}\cdots\alpha_{2}\alpha_{1},$$
where $\alpha_{i}=\beta_s(r_{t-1})$ with probability $p_s$,
and observe that
$$ r_{t}(s^t)=r^{\alpha(t)}.$$
Observe that 
$$\ln\alpha(T)= \sum_{t=1}^{T}\ln\alpha_t.$$
Since 
$$E(\ln\alpha(t)|s^{t-1})-\ln\alpha(t-1)=E(\ln\alpha_t|s^{t-1})\le 0,$$
and 
$$0<b_s \le \beta_s(r_t)\le B_s<\infty .$$
the sequence $Z_T=\ln\alpha(T)$ forms a supermartingale with bounded increments. 
Doob's decomposition theorem gives the representation
$$Z_T=W_T +S_T,$$
where $W_T=\sum_{t=1}^{T}E(\ln\alpha_t|s^{t-1})$ is a decreasing predictable sequence
and $$S_T=\sum_{t=1}^{T} [\ln\alpha_t-E(\ln\alpha_t|s^{t-1})],$$
is a 0 mean martingale with bounded increments.
By Theorem  5.1 (Ch. VII) of \cite{Sh} with probability 1 process $S_T$ either converges to finite limit
or $\limsup_{T\to\infty} S_T=-\liminf_{T\to\infty} S_T=\infty$. In the first case the process $Z_T$ is bounded from above. We will consider only the second case to show  that with positive probability the process $Z_T$ is bounded from above for a set of indices $T$
which has positive density in $\mathbb N$, i.e, there exist  $M>0$, $0<a,b<1$ such that
\begin{equation}\label{Zt}
\pi_{\bf p}(\limsup_{T\to\infty}\frac { \#\{t\le T: Z_t\le M\}}T\ge a)>b.
\end{equation}
Let us denote $$X_t=\ln\alpha_t-E(\ln\alpha_t|s^{t-1})\ ,\ \ t\ge 1.$$ 
This sequence satisfies assumptions of Theorem \ref{Th:Drogin}, with
$\mathcal A_t=\sigma(s^t)$. We have 
\begin{equation*}\begin{split}E(X_t^2|s^{t-1})&=E((\ln\alpha_t-E(\ln\alpha_t|s^{t-1}))^2|s^{t-1})\\
&=
\sum_{s=1}^L p_s(\ln \beta_s(r))^2-\left(\sum_{s=1}^L p_s\ln \beta_s(r)\right)^2={\rm var}(\ln\alpha_t|s^{t-1}) \ge d>0.
\end{split}
\end{equation*}
Thus, the sequence $X_t$ satisfies assumptions of Proposition \ref{arcsine}. In particular, (\ref{whatweneed}) holds,
 i.e., if ${\rm Pos}_T$ is the number of times
$\ln\alpha(t)> 0$ for $t \le T$, then 
$$\limsup_{T \to\infty} [\pi_{\bf p}(\frac {{\rm Pos}_T}T\le a)]=b>0,$$
where $a,b$ are some numbers in $(0,1)$.
This means that if ${N}_T$ is the number of times $\ln\alpha(t) \le 0$ for $t\le T$, then 
$$\limsup_{T\to\infty}[ \pi_{\bf p}(\frac {{N}_T}T \ge 1-a)]=b>0.$$
Now, we  we show that
$$\pi_{\bf p}(\limsup_{T\to\infty}\frac {{N}_T}T\ge 1-a)\ge b/2>0.$$
For $T>T_0$ we have $\pi_{\bf p}(\frac {{N}_T}T \ge 1-a)>b/2$.
Let $A_T=\{\frac {{N}_T}T \ge 1- a\}$, $T\ge T_0$. 
The set which contains points from infinitely many $A_T$ is $A=\cap_{i}\cup_{T>i}A_T$ and
since the sequence $(\cup_{T>i}A_T)_i$ is decreasing we have 
$$\pi_{\bf p}(A)=\lim_{i\to\infty}\pi_{\bf p} (\cup_{T>i}A_T)\ge b/2\ .$$
Thus, with a positive probability $b/2>0$, there exist a sequence
$T_n\to\infty$ such that $\frac {{N}_{T_n}}{T_n} \ge 1-a$ or 
$$\pi_{\bf p}(\limsup_{T\to\infty}\frac {{N}_T}T\ge 1- a)\ge b/2>0.$$
Thus,
 $\ln\alpha(T)$ is negative with positive density, i.e.,
$$\lim_{T\to\infty}\frac 1T \#\{t\le T: \ln\alpha(t)\le 0\}\ge 1- a>0,$$
with positive probability $b/2$. This implies that $r_T(s_T)\ge \bar r>0$ with positive density $1-a$ and positive probability $b/2$. We can construct a 
continuous function $f$ which is 0 around 0 and 1 above $\bar r$. The averages of this function satisfy
$$\limsup_{T\to\infty}\frac 1 T \sum_{T'\le T} f( r_{T'}(s_{T'}))\ge 1- a,$$
with nonzero probability $b/2$ which proves there is no weak convergence to $\delta_0$.
\end{proof}
\begin{remark}
If $E(\ln\alpha_t|s^{t-1})\ge 0$ and ${\rm var}(\ln\alpha_t|s^{t-1})\ge d > 0$, for all $t\ge 1$, using essentially the proof of Theorem \ref{NotSRB}, we obtain that $\delta_1$ is not an SRB measure of $F$. In particular, if $E(\ln\alpha_t|s^{t-1})= 0$ and ${\rm var}(\ln\alpha_t|s^{t-1})\ge d > 0$, for all $t\ge 1$, we obtain that neither $\delta_0$ nor $\delta_1$ is an SRB measure. 
\end{remark}
\section{Properties of $\delta_0$ and $\delta_1$: The case when ${\bf p}$ is unknown}\label{information}
In general, one cannot decide whether $\delta_0$ or $\delta_1$ is the unique SRB measure without having information about ${\bf p}$. We illustrate this fact in the following example.
\begin{example}\label{Ex1}
Let $F=\{\tau_1,\tau_2; p_1,p_2\}$ where
$\tau_1=r^2$, $\tau_2=\sqrt r$ and $p_1, p_2$ are unknown. Observe that the exponents, which are explicit in this case and independent of $r$, are $\beta_1(r)=2$ and $\beta_2(r)=1/2$. Then
$$p_1\ln B_1(r)+p_2\ln B_2(r)= (2p_1-1)\ln 2.$$
By Corollary \ref{Cor1}, if $p_1<1/2$ the measure $\delta_1$ is the unique SRB measure of $F$; however, if $p_1>1/2$ the measure $\delta_0$ is the unique SRB measure of $F$. Thus, for this example, without having information about ${\bf p}$, no information about the nature of $\delta_0$ or $\delta_1$ can be obtained.
\end{example}

Although Example \ref{Ex1} shows that the analysis cannot be definitive in some cases without knowing the probability distribution on $S$, our aim in this section is to find situations when $\delta_0$ and/or $\delta_1$ are not SRB. Moreover, in addition to studying the properties of $\delta_0$ and $\delta_1$,  we are going to study the case when the IFS admit an invariant probability measure whose support is separated from zero and is not necessarily concentrated at one.  The definition of such a measure is given below.
\begin{definition}
Let $\mu$ be a probability measure on $([0,1],\mathfrak B)$, where $\mathfrak B$ is the Borel $\sigma$-algebra.  We define the support of $\mu$, denoted by supp($\mu$), as the smallest closed set of full $\mu$ measure.  We say that supp($\mu$) is separated from zero if there exists an $\eta>0$ such that $\mu([0,\eta])=0$.
\end{definition}
In addition to properties (A) and (B), we assume in this section that:\\

\noindent (C) Every $\tau_s$ has a finite number of fixed points.\\
 
In this section, we use a graph theoretic techniques to analyze ergodic properties of 
$\delta_0$ and $\delta_1$. This approach is inspired by the concept of a \textit{Markov partition} used in the 
dynamical systems literature. For instance, in \cite{BG}, the ergodic properties of a deterministic system which 
admits a Markov partition is studied via a directed graph and an incidence matrix. In our approach we construct a 
partition for our random dynamical system akin to that of a Markov partition and use two directed graphs to study 
ergodic properties of the system.\\ 

We now introduce the two graphs, $G_d$ and $G_u$, which we will use in our analysis. 
\begin{enumerate}
\item Both $G_d$ and $G_u$ have the same vertices;
\item For $s\in\{1,\dots, L\}$, an interval $J_{s,m}=(a_{s,m},a_{s,m+1})$ is a vertex in $G_d$ and $G_u$ if and only if $\tau_s(a_{s,m})=a_{s,m}$,  
$\tau_s(a_{s,m+1})=a_{s,m+1}$ and $\tau_s(r)\not=r$ for all $r\in (a_{s,m},a_{s,m+1})$;
\item Let $J_{s,m}$ and $J_{l,j}$ be two vertices of $G_d$. There is a directed edge connecting $J_{s,m}$ to 
$J_{l,j}$ if and only if $\exists$ an $r\in J_{s,m}$, $r>a_{l,j+1}$, and a $t\ge 1$ such that $\tau_{s}^{t}(r)\in J_{l,j}$. 
\item Let $J_{s,m}$ and $J_{l,j}$ be two vertices of $G_u$. There is a directed edge  connecting $J_{s,m}$ to $J_{l,j}$ 
if and only if $\exists$ an $r\in J_{s,m}$, $r<a_{l,j}$, and a $t\ge 1$ such that $\tau_{s}^{t}(r)\in J_{l,j}$.    
\item  By the \textit{out-degree} of a vertex we mean the number of outgoing directed edges from this vertex in the graph, 
and by the \textit{in-degree} of a vertex we mean the number of incoming directed edges incident to this vertex in the 
graph. 
\item A vertex is called a \textit{source} if it is a vertex with in-degree equals to zero. 
A vertex is called a \textit{sink} if it is a vertex with out-degree equals to zero.
\end{enumerate}
For the above graphs, one can identify two types of vertices: 
let $(a_{s,m},a_{s,m+1})$ be a vertex. If $\tau_s(r)>r$ for all $r\in(a_{s,m},a_{s,m+1})$, 
then the vertex $(a_{s,m},a_{s,m+1})$ will be denoted by $\hat J_{s,m}$. 
If $\tau_s(r)<r$ for all $r\in(a_{s,m},a_{s,m+1})$, then the vertex $(a_{s,m},a_{s,m+1})$ 
will be denoted by $\check J_{s,m}$. When we prove a statement for a vertex $J_{s,m}$ (without a label), 
this means that the result holds for both types of vertices. The following lemma contains some properties of $G_d$ and 
$G_u$.
\begin{lemma}\label{Le1} Let $G_d$ and $G_u$ be defined as above. 
\begin{enumerate}
\item If $\hat J_{s,m}$ is a vertex in $G_d$, then $\hat J_{s,m}$ is a sink in $G_d$. 
\item Let $\check J_{s,m}$ and $J_{l,j}$ be two vertices in $G_d$. There is a directed edge connecting 
$\check J_{s,m}$ to $J_{l,j}$ in $G_d$ if and only if $a_{s,m}< a_{l,j+1}<a_{s,m+1}$. 
In particular for all $s\in S$ there is no directed edge in $G_d$ connecting $J_{s,m}$ to $J_{s,j}$ for any $m$ and $j$.
\item If $\check J_{s,m}$ is a vertex in $G_u$, then $\check J_{s,m}$ is a sink in $G_u$. 
\item Let $\hat J_{s,m}$ and $J_{l,j}$ be two vertices in $G_u$. 
There is a directed edge connecting  $\hat J_{s,m}$ to $J_{l,j}$ in $G_u$ if and only if 
$a_{s,m}< a_{l,j}<a_{s,m+1}$. In particular for all $s\in S$ there is no directed edge in $G_u$ 
connecting $J_{s,m}$ to $J_{s,j}$ for any $m$ and $j$.
\end{enumerate}
\end{lemma}
\begin{proof}
The proof of the first statement is straight forward. 
Indeed, let $J_{l,j}$ be any vertex in $G_d$ and $r\in\hat J_{s,m}$ such that $r>a_{l,j+1}$. 
Then for all $t\ge 1$ $\tau_s^t(r)>\tau_s^{t-1}(r)>\dots\tau_s(r)>r>a_{l,j+1}$. 
The proof of the second statement follows from the fact that if $r>a_{s,m}\ge a_{l,j+1}$ then for 
$t\ge1$ we have $\tau_s^{t}(r)>a_{s,m}\ge a_{l,j+1}$. If $r>a_{l,j+1}>a_{s,m}$, then there exits a 
$t\ge1$ such that $a_{s,m}<\tau_s^t(r)<a_{l,j+1}$. 
Proofs of the third and fourth statements are similar to the first two.  
\end{proof}
For our further analysis we introduce the following notion.
\begin{definition}
We say that a random orbit of $F$ stays above a point $c$ if all the points of the infinite orbit are bigger than 
or equal to $c$ with probability $\pi_{\bf p}=1$. 
\end{definition}
\begin{lemma}\label{Le2}
Let $J_{l,j}$ be a vertex in $G_d$ such that $a_{l,j+1}\not=1$. 
If $J_{l,j}$ is a source in $G_d$, then the random orbit of $F$ starting from $r>a_{l,j+1}$ stays above 
$a_{l,j+1}$ with probability $\pi_{\bf p}=1$.
\end{lemma}
\begin{proof}
Suppose $J_{l,j}$ is a source in $G_d$. Then for all 
$r>a_{l,j+1}$, we have $\tau^t_s(r)>a_{l,j+1}$ for all $s\in S$ and $t\ge 1$. This means that if $r>a_{l,j+1}$ we have 
$\tau_{s_1}(r)>a_{l,j+1}$ and $\tau_{s_2}\circ\tau_{s_1}(r)>a_{l,j+1}$ and so on. 
\end{proof}
\begin{theorem}\label{prop1} Let $F$ be an IFS whose transformations satisfy the properties (A), (B) and (C).
\begin{enumerate}
\item If for all $s\in S$ there is a vertex $\check J_{s,m}$ in $G_d$ with $a_{s,m}=0$, 
then $\delta_0$ is an SRB measure, $B(\delta_0)\supseteq[0,a)$, where $a=\min_s\{a_{s,m+1}\}$. 
In particular, for any $r_0\in [0,a)$, $\lim_{t}r_t(s^t)=0$ a.s.
\item If for all $s\in S$ there is a vertex $\hat J_{s,m}$ in $G_d$ with $a_{s,m+1}=1$, 
then $\delta_1$ is an SRB measure. Moreover, $B(\delta_1)\supseteq (b,1]$, 
where $b=\max_s\{a_{s,m}\}$. In particular, for any $r_0\in (b,1]$, $\lim_{t}r_t(s^t)=1$ a.s.
\item Let $J_{l,j}$ be a vertex in $G_d$ such that $a_{l,j+1}\not=1$. 
If $J_{l,j}$ is a source in $G_d$\footnote{In the case where $a_{l,j}=0$, even if other $\hat J_{s,m}$, with $a_{s,m}=0$, 
receives a directed edge, the result still holds. Thus, to know the existence of an invariant probability measure whose 
support is separated from 0, it is enough to check that one vertex $J_{l,j}$ with $a_{l,j}=0$ which is a source in $G_d$. 
Statements of Lemma \ref{Le1} can be useful to visualize cases of this type.} then $F$ preserves a probability measure 
whose support is separated from 0 \footnote{The invariant measure here is not necessarily $\delta_1$.}. 
\item Let $J_{l,j}$ be a vertex in $G_u$ such that $a_{l,j+1}\not=0$. If $J_{l,j}$ is a source in $G_u$ then $F$ 
preserves a probability measure whose support is separated from 1. 
\item Let $\hat J_{s^*,m}$ be a vertex with $a_{s^*,m}=0$ whose out-degree in $G_u$ is at least one. If $\exists$ 
a vertex $J_{s_0,j}$ in $G_d$, $a_{s_0,j}=0$ and $a_{s_0,j+1}<a_{s^*,m+1}$, which is a source in $G_d$, then for any 
$r_0\in (0,1]$, $\lim_{t}r_t(s^t)\not=0$ a.s. Moreover, $\delta_0$ is not an SRB measure for $F$.
\item Let $\check J_{s_0,j}$ be a vertex in $G_d$ such that $a_{s_0,j+1}=1$ and whose out-degree in $G_d$ is at least one. 
If $\exists$ a $J_{s^*,m}$ in $G_u$, $a_{s^*,m+1}=1$ and $a_{s^*,m}>a_{s_0,j}$, which is a source in $G_u$, 
then for any $r_0\in [0,1)$, $\lim_{t}r_t(s^t)\not=1$ a.s. Moreover, $\delta_1$ is not an SRB measure for $F$.
\item If for all $s\in S$ the vertices whose $a_{s,m}=0$ are of the form $\hat J_{s,m}$ and their $a_{s,m+1}\equiv a$ 
are identical, then for any $r_0$ in $(0,a]$, with probability one, $\lim_{t}r_t(s^t)=a$. In particular, $\delta_a$ 
is an SRB measure with $B(\delta_a)=(0,a]$ and $\delta_0$ is not an SRB measure. 
\item If for all $s\in S$ the vertices whose $a_{s,m+1}=1$ are of the form $\check J_{s,m}$ and their $a_{s,m}\equiv b$ 
are identical, then for any $r_0$ in $[b,1)$, with probability one, $\lim_{t}r_t(s^t)=b$. In particular, $\delta_b$ is an 
SRB measure with $B(\delta_b)=[b,1)$ and $\delta_1$ is not an SRB measure. 
\end{enumerate}
\end{theorem}    
\begin{proof}
We only prove the odd numbered statements in the theorem. Proofs of the even numbered statements are very similar. \\
(1) For any $r_0\in[0,a)$, any random orbit of $F$ will converge to zero. Using Lemma \ref{fix}, this shows that $\delta_0$ is an SRB measure with $B(\delta_0)\supseteq[0,a)$. \\  
(3) Let $r_0>a_{l,j+1}$.  Since $[0,1]$ is a compact metric space and for all $s\in S$ $\tau_s$ is continuous, the average $\frac{1}{T}\sum_{t=0}^{T-1}P^t\delta_{r_0}$ of the probability measures converges in the weak* topology to an $F$ invariant probability measure\footnote{This follows from a random version of the Krylov-Bogoliubov Theorem \cite{Ar}.}. By Lemma \ref{Le2}, this measure is supported on $[a_{l,j+1},1]$. \\ 
(5) Let $D=\{J_{s,m}\setminus\{0\}:\, a_{s,m}=0\}$. For any $r_0\in D$, there exists a finite $t\ge 1$ such 
that $\tau_{s^*}^t(r_0)>a_{s_0,j+1}$. Since $J_{s_0,j}$ is a source in $G_d$, by Lemma \ref{Le2}, $\tau_{s^*}^t(r_0)$ stays above $a_{s_0,j+1}$ with probability $\pi_{\bf p}=1$. Therefore, for any $r_0\in D$,  with positive probability, the random orbit of $r_0$ is bounded away from $0$. Let us consider now the case of a starting point $r_0'>a_{s_0,m+1}$. Since all the transformations are homeomorphisms and $0$ is a common fixed point, for any $r'_0> a_{s_0,m+1}$ and any $t\ge 0$, with positive probability, $r_t(s^t)>a_{s_0,m+1}$. Hence, for any $r\in (0,1]$, with strictly positive probability, $\lim_{t\to\infty}r_t(s^t)\ge a_{s_0,m+1}$. Moreover, with strictly positive probability, for any $r\in (0,1]$, there exists a $T-1>t_0\ge1$ such that
$$\frac{1}{T}\sum_{t=0}^{T-1}r_{t}(s^t)\ge \frac{1}{T}\sum_{t=0}^{t_0-1}r_{t}(s^t)-\frac{(t_0+1)}{T}a_{s_0,m+1}+a_{s_0,m+1}.$$
Therefore, with strictly positive probability, for any $r\in(0,1]$,
\begin{equation}\label{NSRB}
\lim_{T\to\infty}\frac{1}{T}\sum_{t=0}^{T-1}r_{t}(s^t)\ge a_{s_0,m+1}.
\end{equation}
Now, to show that $\delta_0$ is not an SRB measure, it is enough to find a continuous function $f$ on $[0,1]$ such that with positive probability, for any $r\in(0,1]$,
$$\{\lim_{T\to\infty}\frac{1}{T}\sum_{t=0}^{T-1}f(r_t(s^t))\not=f(0)\}.$$
Indeed, this is the case if we use the function $f(r)=r$ and (\ref{NSRB}). Thus, $\delta_0$ is not an SRB measure.. \\
(7) Obviously, for any $r_0\in(0,a]$, the random orbit of $F$ starting at $r_0$ will converge to $a$. Using Lemma \ref{fix}, this implies that $\delta_a$ is an SRB measure with $B(\delta_a)=(0,a]$. Moreover, since all the transformations are homeomorphisms with common fixed point at $a$, for any $r_0'>a$, the random orbit of $F$ stays above $a$.  Thus, $\delta_0$ is not an $SRB$ measure. 
\end{proof}
\section{Asset Market Games}\label{games}
In this section, we apply our results to \textit{evolutionary models of financial markets}. In particular, we will focus on the model introduced by \cite{EHS}. First, we recall the model of \cite{EHS}.

\subsection{The Model} Let $S$ is a finite set and $s_{t}\in S$, $t=1,2,...$, be the ``state of the world'' at
date $t$. Let ${\bf p}$ be a probability distribution on $S$ such that for all $s\in S$ ${\bf p}(s)>0$. We also assume that $s_{t}$ are independent and identically distributed. 

\bigskip 

In this model there are $K$ {\it ``short-lived'' assets} $k=1,2,...,K$ (live one period and are identically reborn every next period). One unit of asset $k$ issued at time $t$ yields payoff $D_{k}(s_{t+1})\geq 0$ at time $t+1$. It is assumed that 
\[
\sum\nolimits_{k=1}^{K}D_{k}(s)>0\ \ \text{for all \ }s\in S 
\]
and 
\[
ED_{k}(s_{t})>0 
\]
for each $k=1,2,...,K$ , where $E$ is the expectation with respect to the
underlying probability ${\bf p}$. The total amount of
asset $k$ available in the market is $V^{k}=1.$

\bigskip 

In this model there are $I$ {\it investors} {\it (traders)} $i=1,...,I$. Every investor $i$ at each time $t=0,1,2,...$ has a {\it portfolio} 
\[
x_{t}^{i}=(x_{t,1}^{i},...,x_{t,K}^{i}), 
\]
where $x_{t,k}^{i}$ is the number of units of asset $k$ in the portfolio $x_{t}^{i}=x_{t}^{i}(s^{t}),\;s^{t}=(s_{1},...,s_{t})$. We assume that for each moment of time $t\geq 1$ and each random situation $s^{t}$, the market for every asset $k$ clears: 
\begin{equation}
\sum_{i=1}^{I}x_{t,k}^{i}(s^{t})=1. 
\end{equation}
Each investor is endowed with initial wealth $w_{0}^{i}>0.$ Wealth $w_{t+1}^{i}$ of investor $i$ at time $t+1$
can be computed as follows: 
\begin{equation}\label{wealth1}
w_{t+1}^{i}=\sum_{k=1}^{K}D_{k}(s_{t+1})x_{t,k}^{i}. 
\end{equation}
{\it Total\ market wealth }at time $t+1$ is equal to 
\begin{equation}\label{total}
w_{t+1}=\sum_{i=1}^{I}w_{t+1}^{i}=\sum_{k=1}^{K}D_{k}(s_{t+1}).
\end{equation}
Investment strategies are characterized in terms of {\it investment proportions:}
$$\Lambda ^{i}=\{\lambda _{0}^{i},\,\lambda _{1}^{i},\,\lambda _{2}^{i},...\}$$
of $K$-dimensional vector functions $\lambda _{t}^{i}=(\lambda _{t,1}^{i},...,\lambda _{t,K}^{i}),\;\lambda
_{t,k}^{i}=\lambda _{t,k}^{i}(s^{t})\;t\geq 0,$ satisfying $\lambda _{t,k}^{i}>0,\;\sum_{k=1}^{K}\lambda _{t,k}^{i}=1.$ Here, $\lambda _{t,k}^{i}$ stands for the {\it share of the budget} $w_{t}^{i}$ of investor $i$ that is invested into asset $k$
at time $t$. In general $\lambda _{t,k}^{i}$ may depend on $s^{t}=(s_{1},s_{2},...,s_{t})$. Given strategies $\Lambda ^{i}=\{\lambda _{0}^{i},\;\lambda _{1}^{i},\;\lambda
_{2}^{i},\;...\}$ of investors $i=1,...,I$, the equation 
\begin{equation}
p_{t,k}\cdot 1=\sum_{i=1}^{I}\lambda _{t,k}^{i}w_{t}^{i}\;  
\end{equation}
determines the market clearing {\it price} $p_{t,k}=p_{t,k}(s^{t})$ of asset $k$. The number of units of asset $k$ in the portfolio of investor $i$ at time $t$ is equal to 
\begin{equation}
x_{t,k}^{i}=\frac{\lambda _{t,k}^{i}w_{t}^{i}}{p_{t}^{k}}.  
\end{equation}
Therefore 
\begin{equation}\label{port}
x_{t,k}^{i}=\frac{\lambda _{t,k}^{i}w_{t}^{i}}{\sum_{j=1}^{I}\lambda
_{t,k}^{j}w_{t}^{j}}.  
\end{equation}
By using (\ref{port}) and (\ref{wealth1}), we get 
\begin{equation}\label{wealth2}
w_{t+1}^{i}=\sum_{k=1}^{K}D_{k}(s_{t+1})\frac{\lambda _{t,k}^{i}w_{t}^{i}}{
\sum_{j=1}^{I}\lambda _{t,k}^{j}w_{t}^{j}}. 
\end{equation}
Since $w_{0}^{i}>0$, we obtain $w_{t}^{i}>0$ for each $t$. The main focus of the model is on the analysis of the
dynamics of the {\it market shares} of the investors 
\[
r_{t}^{i}=\frac{w_{t}^{i}}{w_{t}},\;i=1,2,...,I. 
\]
Using (\ref{wealth2}) and (\ref{total}), we obtain 
\begin{equation}\label{relative}
r_{t+1}^{i}=\sum_{k=1}^{K}R_{k}(s_{t+1})\frac{\lambda
_{t,k}^{i}r_{t}^{i}}{\sum_{j=1}^{I}\lambda _{t,k}^{j}r_{t}^{j}},\;i=1,2,...,I, 
\end{equation}
where 
\[
R_{k}(s_{t+1})=\frac{D_{k}(s_{t+1})}{\sum_{m=1}^{K}D_{m}(s_{t+1})} 
\]
are the {\it relative }({\it normalized}){\it \ payoffs} of the assets $
k=1,2,...,K$. We have $R_{k}(s)\geq 0$ and $\sum_{k}R_{k}(s)=1$.

\subsection{Performance of investment strategies and the Kelly rule} In the theory of \textit{evolutionary finance} there are three possible grades for investor $i$ (or for the
strategy she/he employs):

(i) {\it extinction}: $\lim r_{t}^{i}=0\, \text{a.s.}$;

(ii) {\it survival}: $\lim \sup r_{t}^{i}>0$ but $\lim \inf r_{t}^{i}<1$\, \text{a.s.};

(iii) {\it domination}: $\lim r_{t}^{i}=1$\, \text{a.s.}

\begin{definition}
An investment strategy is called \textit{completely mixed strategy} if it assigns a positive percentage of wealth 
$\lambda _{t,k}(s^{t})$ to every asset $k=1,\ldots ,K$ for all $\ t$ and $s^{t};$ moreover, it is called \textit{simple }
if $\lambda_{t,k}(s^{t})=\lambda _{k}>0.$
\end{definition}
In this theory, the following \textit{simple} portfolio rule has been very successful: define 
\[
\lambda ^{\ast }=(\lambda _{1}^{\ast },...,\lambda _{K}^{\ast }),\;\lambda
_{k}^{\ast }=ER_{k}(s_{t}),\;k=1,...,K, 
\]
so that $\lambda _{1}^{\ast },...,\lambda _{K}^{\ast }$ are the {\it 
expected relative payoffs} of assets $k=1,...,K$. The portfolio rule $\lambda ^{\ast }$ is called the \textit{Kelly rule} which expresses the
investment principle of ``betting your beliefs'' \cite{Ke}. In \cite{EHS} under the following two conditions:\\

\noindent E1) There are no {\it redundant assets}, i.e. the functions $
R_{1}(s),...,R_{K}(s)$\ of $s\in S$\ are linearly independent.\\

\noindent E2) All investors use {\it simple strategies;}\\

\noindent it was shown that investors who follow the Kelly rule survive and 
others who use a different simple strategy get extinct. In particular, If only one investor follows the Kelly rule,
then this investor dominates the market.

\bigskip

The main challenge in using the Kelly rule lies in the fact that it requires from investors the full knowledge of the probability distribution ${\bf p}$. In Subsection \ref{pinformation}, using an IFS representation of (\ref{relative}) and Theorem \ref{main}, we overcome this difficulty by finding another successful strategy which requires partial knowledge of the probability distribution ${\bf p}$.
\subsection{An IFS realization of the model}
In the rest of the paper, we are going to show how the above model can be represented by an IFS. We are going to apply the results of Sections \ref{martingale} and \ref{information} to study the dynamics of (\ref{relative}). As in \cite{EHS}, we assume here that all the investors use simple strategies. Further, we focus on the case\footnote{This is the same as assuming that there are $I$ investors, $I>2$, where $I-1$ investors use the same strategy and only one investor deviates from them.} when $I=2$. The market selection process (\ref{relative}) reduces to the following one dimensional system:
\begin{equation}\label{reduced}
r_{t+1}(s^{t+1})=\sum_{k=1}^{K}R_{k}(s_{t+1})\frac{\lambda _k^1 r_t}{\lambda _k^1 r+\lambda _k^2(1- r_t)},
\end{equation}
where $r_t$ is investor's 1 relative market share at time $t$ and $(\lambda_k^1)_{k=1}^{K}$ and $(\lambda_k^2)_{k=1}^{K}$ are the investment strategies of investor 1 and 2 respectively. Then the random dynamical (\ref{reduced}) of the market selection process can be described by an iterated function system with probabilities:
$$F=\{\tau_1,\tau_2,\dots ,\tau_L;p_1,p_2,\dots ,p_L\},$$
where
$$\tau_s(r)=\sum_{k=1}^{K}R_{k}(s)\frac{\lambda _k^1 r}{\lambda _k^1 r+\lambda _k^2(1- r)}.$$ 
We first note that the transformations $\tau_s$ of the IFS of the market selection process are maps 
from the unit interval into itself and they satisfy assumptions (A), (B) and (C). In fact, the maps for this model have additional properties. For example, they are differentiable functions. 
\subsection{Investors with partial information on ${\bf p}$ and a generalization of the Kelly rule}\label{pinformation}
We use Theorem \ref{main} to provide a rule for investors with partial information on ${\bf p}$. The investor who 
follows this rule cannot be driven out of the market; i.e., she/he either \textit{dominates} or at least \textit{survives}. The importance of this rule lies in the fact that investor 1 does not need to know the Kelly rule exactly\footnote{It is often difficult for an investor to know the exact probability distribution of the states of the world.}. She/he only needs to know a perturbation of the Kelly rule; for example, the Kelly rule plus some error bounds. 

Firstly, we show in the following lemma that the logarithms of the exponents $\beta_{s}(r)$ are uniformly bounded.
\begin{lemma}\label{es} Let 
$$\tau(r)=\sum_{k=1}^{K}R_{k}\frac{\lambda _k^1 r}{\lambda _k^1 r+\lambda _k^2(1- r)}\ , \ \  r\in [0,1] ,$$ 
and $$\tau(r)=r^{\beta(r)},$$
where, for each $1\le k\le K$ we have $R_k\ge 0$, $\lambda_k^1> 0$, $\lambda_k^2> 0$ and
 $\sum_{k=1}^K R_k=\sum_{k=1}^K \lambda^1_k=\sum_{k=1}^K \lambda^2_k=1$.
Then for any $r\in U$, $U\subseteq [0,1]$, $\ln(\beta(r))$ is bounded. 
\end{lemma}
\begin{proof}
Without loss of generality, we assume that $U=[0,1]$. We have $\tau(r)=r^{\beta(r)}= \exp(\ln(r)\beta(r))$, so 
$$
\beta(r)=\frac {\ln(\tau(r))}{\ln(r)}.
$$
The minimum and maximum of $\beta(r)$ can be attained at $r=0$, $r=1$ or at a point of a local extremum.
Using De L'Hospital rule we find 
$$\lim _{r\to 0^+}\beta(r)=1 \ \ \ \ \ \mathrm{and} \ \ \ \ \ \lim_{r\to 1^-}\beta(r)=\sum_{k=1}^{K}R_k\frac{\lambda_k^2}{\lambda_k^1}.$$
A point of local extremum $r_{\ast}$ in $(0,1)$ of $\beta(r)$ is found by solving
$$\beta'(r)=\frac{1}{\ln(r)} \left(\frac{\tau'(r)}{\tau(r)}-\frac{B(r)}{r}\right)=0.$$
Therefore, at the point $r=r_{\ast}$ of local extremum
$$\beta(r_\ast)=\frac{\sum_{k=1}^{K}R_k\frac{\lambda_k^1\lambda_k^2}{[\lambda_k^1r_\ast+\lambda_k^2(1-r_\ast)]^2}}{\sum_{k=1}^{K}R_k\frac{\lambda_k^1}{\lambda_k^1r_\ast+\lambda_k^2(1-r_\ast)}}.$$
Observe that the function 
$$\frac{\sum_{k=1}^{K}R_k\frac{\lambda_k^1\lambda_k^2}{[\lambda_k^1r+\lambda_k^2(1-r)]^2}}{\sum_{k=1}^{K}R_k\frac{\lambda_k^1}{\lambda_k^1r+\lambda_k^2(1-r)}}$$
is continuous at $[0,1]$. Thus, it attains its maximum and minimum on $[0,1]$. This completes the proof of the lemma.
\end{proof}
\begin{corollary}
Let 
$$\tau_s(r)=\sum_{k=1}^{K}R_{k}(s)\frac{\lambda _k^1 r}{\lambda _k^1 r+\lambda _k^2(1- r)}\ , \ \  r\in [0,1] .$$
Then for $r\in U$, $U\subseteq [0,1]$,
$$b_s=\min_{r\in \bar{U}}\frac{\sum_{k=1}^{K}R_k\frac{\lambda_k^1\lambda_k^2}{[\lambda_k^1r+\lambda_k^2(1-r)]^2}}{\sum_{k=1}^{K}R_k\frac{\lambda_k^1}{\lambda_k^1r+\lambda_k^2(1-r)}}\ \ \ \ \ \mathrm{and} \ \ \ \ \ B_s=\max_{r\in \bar{U}}\frac{\sum_{k=1}^{K}R_k\frac{\lambda_k^1\lambda_k^2}{[\lambda_k^1r+\lambda_k^2(1-r)]^2}}{\sum_{k=1}^{K}R_k\frac{\lambda_k^1}{\lambda_k^1r+\lambda_k^2(1-r)}}.$$
\end{corollary} 
\begin{theorem}
If for each $k\in\{1,\dots,K\}$ $\lambda^1_k$ lies between $ER_k$ and $\lambda_k^2$, then investor 1 cannot be driven out of the market; i.e., she/he either \textit{dominates} or at least \textit{survives}. 
\end{theorem}
\begin{proof}
Let us consider the function
$$G(r)=\sum_{k=1}^K  v_k \frac {\lambda_k^1}{\Lambda_k(r)}\ \ , \ \ r\in[0,1],$$
where $$\Lambda_k(r)=\lambda_k^1 r +\lambda_k^2(1-r)=(\lambda_k^1-\lambda_k^2)r+\lambda_k^2,$$
and $V=(v_1,v_2,\dots,v_L)$ is a probability vector. We will find conditions on $\lambda^1_k$ which ensure $G(r)\ge 1$,
$r\in[0,1]$. It is easy to see that
\begin{equation}\label{one}
G(1)=\sum_{k=1}^K  v_k \frac {\lambda_k^1}{\lambda_k^1}=1.
\end{equation}
We also have
$$G'(r)=\sum_{k=1}^K  v_k {\lambda_k^1}\frac {-(\lambda_k^1-\lambda_k^2)}{(\Lambda_k(r))^2}\ \ , \ \ r\in[0,1],$$
and
$$G''(r)=\sum_{k=1}^K  v_k {\lambda_k^1}\frac {2(\lambda_k^1-\lambda_k^2)^2}{(\Lambda_k(r))^3}>0\ \ , \ \ r\in[0,1].$$
Thus, $G$ is a convex function and its derivative $G'$ is increasing. If $G'(1)\le 0$ then $G$ is decreasing and
because of (\ref{one}) this implies that $G(r)\ge 1$, $r\in[0,1]$. Observe that
$$G'(1)=\sum_{k=1}^K    v_k {\lambda_k^1}\frac {-(\lambda_k^1-\lambda_k^2)}{(\lambda_k^1)^2}=
\sum_{k=1}^K \frac{v_k}{\lambda_k^1}  {(\lambda_k^2-\lambda_k^1)} .$$
It is easy to see that a sufficient condition for $G'(1)\le 0$ is
\begin{equation}\label{generalized}\begin{split}
v_k\ge \lambda_k^1   \ \  \ \text {   if   }\ \ \ \lambda_k^1 \ge \lambda_k^2;\\ 
v_k\le \lambda_k^1   \ \  \ \text {   if   }\ \ \ \lambda_k^1 \le \lambda_k^2,
\end{split}
\end{equation}
or, in short, for each $k$, $1\le k\le K$, $ \lambda_k^1$ should be between $\lambda_k^2$ and $v_k$.

Now, let us consider the expression 
$$\sum_{s=1}^L p_s \ln(\beta_s(r)).$$ 
We have
\begin{equation}\begin{split}
&\sum_{s=1}^L p_s \ln(\beta_s(r))\le \ln\left(\sum_{s=1}^L p_s \beta_s(r)\right)
=\ln\left(\sum_{s=1}^L p_s \frac{\ln(\tau_s(r))}{\ln r}\right)\\
&\le\ln\left(\frac 1{\ln r}\ln\left( \sum_{s=1}^L p_s \tau_s(r)\right)\right)=
\ln\left(\frac 1{\ln r}\ln\left( \sum_{s=1}^L p_s \sum_{k=1}^K  R_k(s) \frac {\lambda_k^1r}{\Lambda_k(r)}\right)\right)\\
&=\ln\left(\frac 1{\ln r}[\ln r+\ln\left(  \sum_{k=1}^K  (\sum_{s=1}^L p_s R_k(s)) \frac {\lambda_k^1}{\Lambda_k(r)}\right)]\right)
=\ln\left(1+\frac 1{\ln r}\ln( G(r))\right),
\end{split}
\end{equation}
with $v_k$ being the expected payoff for the $k^{\text{th}}$ asset, $v_k=\sum_{s=1}^L p_s R_k(s)$, $k=1,\dots,K$.

A sufficient condition for $\sum_{s=1}^L p_s \ln(\beta_s(r))\le 0$ is  for $r\in [0,1]$:
$$\ln( G(r))\ge 0 \ \ \text{or equivalently}  \ \  G(r)\ge 1. $$
We have shown before that a sufficient condition for this is (\ref{generalized}) or placing each $\lambda_k^1$ between
the expected payoff $v_k$ and $\lambda_k^2$. 

To complete the proof of the theorem, we first use Lemma \ref{es} to observe that exponents $\beta_s(r)$ of this system are 
bounded and then (1) of Theorem \ref{main}. Indeed, for any fixed partial history $s^{t-2}$, because the stochastic process $s_t$ is an iid process, 
we have
$$E(\ln\alpha_{t}|s^{t-1})=\sum_{s=1}^L p_s \ln(\beta_s(r_{t-2})).$$
\end{proof}
\subsection{Incorrect beliefs}
Our results in Section \ref{information} are also interesting for studying the dynamics of (\ref{relative}). In fact, they can be used to study the dynamics in the situation where both players have `incorrect beliefs'; i.e., when players do not have the right information or partial information about ${\bf p}$.  Thus, they either use wrong distributions to build their strategies or they arbitrarily choose their strategies. Consequently, their strategies are, in general, different from the Kelly rule and the generalization which we presented in Subsection \ref{pinformation}. In this case, the results of Section \ref{information} can be used to identify the exact outcome of the game in certain situations. In some situations, as in Example \ref{Ex1}, one cannot know the outcome of the system without knowing ${\bf p}$. 
\section{Appendix}\label{app}
The following general arcsine law has been proved in \cite{Dro}. 
\begin{theorem}\cite{Dro}\label{Th:Drogin}
Let $X_1,X_2,\dots$ be a sequence of random variables adapted to the sequence of $\sigma$-algebras
$\mathcal A_1, \mathcal A_2,\dots$. Let $S_m=\sum_{i=1}^m X_i$, $v_m=\sum_{i=1}^m E(X_i^2|\mathcal A_{i-1})$ and assume
$$ E(X_{m+1}|\mathcal A_m)=0\ ,\ \ \ \ \ EX_m^2<\infty\ ,\ \  \ \text{and}\ \ \ v_m\to \infty \ \text{a.s.}$$
Let $T_n=\inf\{m:v_m\ge n\}$ and $L_n=\frac 1 n \sum_{i=1}^{T_n} E(X_i^2|\mathcal A_{i-1})\chi_{\{S_i>0\}}$. 
If $$\frac 1n \sum _{i=1}^{T_n}X_i^2 \chi_{\{X_i^2>n\eps\}}\underset{L_1}{\longrightarrow} 0 \ \ \ \text{for all} \ \eps>0,$$
then the distributions of $L_n$ converge to the arcsine distribution.
\end{theorem}

We now use Theorem \ref{Th:Drogin} to prove a proposition which is used in the proof of Theorem \ref{NotSRB}.
\begin{proposition}\label{arcsine} Let $X_1,X_2,\dots$ be a sequence of random variables adapted to the sequence of $\sigma$-algebras
$\mathcal A_1, \mathcal A_2,\dots$. Suppose that there exist constants $d>0$ and $0<D<\infty$ such that for all $n\ge 1$ we have 
$$0<d\le E(X_n^2|\mathcal A_{n-1})\ \  \text{and}\ \ X_n^2\le D.$$ 
Then, the sequence satisfies the remaining assumptions of Theorem \ref{Th:Drogin}. In particular, Theorem \ref{Th:Drogin} implies the condition 
\begin{equation}\label{whatweneed}
\limsup_{n \to\infty} Pr(\frac {{\rm Pos}_n}n\le a)\ge b>0,
\end{equation}
for some constants $0<a,b<1$, where ${\rm Pos}_n=\sum_{i=1}^{n} \chi_{\{S_i>0\}}$.
\end{proposition}
\begin{proof}
The remaining assumptions of 
Theorem \ref{Th:Drogin} are trivially satisfied. We have $m\cdot d\le v_m\le m\cdot D$ for all $m\ge 1$ so 
$T_n\cdot d \le n\le T_n\cdot D$ for all $n\ge 1$. Then,
$$L_n \ge \frac 1 D \frac 1{T_n}  d \sum_{i=1}^{T_n} \chi_{\{S_i>0\}}\ $$
and, for $0\le a_1\le 1$, we have
$$ Pr\left( \frac d D \frac 1{T_n}  \sum_{i=1}^{T_n} \chi_{\{S_i>0\}}\le a_1\right)\ge Pr( L_n\le a_1) \underset{n\to \infty}{\longrightarrow}
\frac 2{\pi}\arcsin \sqrt{a_1}.$$
For $a_1$ small enough we obtain a meaningful estimate
$$ Pr\left(  \frac 1{T_n}  \sum_{i=1}^{T_n} \chi_{\{S_i>0\}}\le a\right)\ge  \frac 1{\pi}\arcsin \sqrt{a},$$
for $a=a_1\frac D d$  and $n$ large enough. This implies condition (\ref{whatweneed}).
\end{proof}


\begin{thebibliography}{99}
\bibitem{Ar} (1723992)
\newblock L. Arnold, 
\newblock ``Random Dynamical Systems," 
\newblock Springer Verlag, Berlin, 1998.
\bibitem {BG} (1461536)
\newblock A. Boyarsky and P. G\'ora, 
\newblock ``Laws of Chaos,'' 
\newblock Brikha\"user, Boston, 1997.
\bibitem{Bu} (1707698)  
\newblock J. Buzzi, 
\newblock \emph{ Absolutely continuous S.R.B. measures for random Lasota-Yorke maps}, 
\newblock Trans. Amer. Math. Soc., {\bf 352} (2000), 3289--3303.
\bibitem{EHS} (1926235) 
\newblock I. Evstigneev, T. Hens and K.R. Schenk-Hopp\'{e}, 
\newblock \emph{Market selection of financial trading strategies: Global stability}, 
\newblock Math. Finance, {\bf 12} (2002), 329--339.
\bibitem{Di} (1669737) 
\newblock P. Diaconis and D. Freedman,
\newblock \emph{Iterated random functions,}
\newblock SIAM Rev., {\bf 41} (1999), 45--76. 
\bibitem{Dro} (0303595)
\newblock R. Drogin,
\newblock \emph{An invariance principle for martingales}, 
\newblock Ann. Math. Statist., {\bf 43} (1972), 602--620.
\bibitem{DF} (0193668)
\newblock L. Dubins and D. Freedman, 
\newblock \emph{Invariant probabilities for certain Markov processes}, 
\newblock Ann. Math. Statist., {\bf 37} (1966), 837--848.
\bibitem{HH} (0624435)
\newblock P. Hall and C. Heyde, 
\newblock ``Martingale Limit Theory and Its Application," 
\newblock Academic Press, New York-London, 1980.  
\bibitem{Ke} (0090494)
\newblock J.L. Kelly, 
\newblock \emph{A new interpretation of information rate}, 
\newblock Bell Sys. Tech. J., {\bf 35} (1956), 917--926.
\bibitem{Ki} (0874051)
\newblock Y. Kifer, 
\newblock ``Ergodic Theory of Random Transformations," 
\newblock  Birkh\"auser, Boston, 1986.
\bibitem{Liu} (1855833)
\newblock P-D. Liu, 
\newblock\emph{Dynamics of random transformations: smooth ergodic theory}, 
\newblock Ergodic Theory Dynam. Syst., {\bf 21} (2001), 1279--1319. 
\bibitem{Sh} (0737192) 
\newblock A.N. Shiryaev, 
\newblock ``Probability,"
\newblock Springer-Verlag, New York, 1984.
\bibitem{St} (1962693) 
\newblock \"O. Stenflo, 
\newblock \emph{Uniqueness of invariant measures for place-dependent random iterations of functions},
\newblock in ``Fractals in Multimedia" (eds. M.F. Barnsley, D. Saupe and E.R. Vrscay), Springer, (2002), 13--32.
\bibitem{Y} (1933431)
\newblock L-S. Young, 
\newblock \emph{What are SRB measures, and which dynamical systems have them?,}  
\newblock J. Statist. Phys., {\bf 108} (2002), 733--754.
\end{thebibliography}
\end{document}